\documentclass[a4paper,12pt,reqno]{amsart}
\usepackage{graphics}
\newtheorem{lem}{Lemma}

\newtheorem{theo}{Theorem}

\numberwithin{equation}{section}

\newcommand{\sgn}{\operatorname{sgn}}

\newcommand{\id}{\operatorname{id}}

\begin{document}

\title[Alternating sign matrices]{A new proof of the refined alternating sign matrix theorem}

\author[Ilse Fischer]{\box\Adr}

\newbox\Adr
\setbox\Adr\vbox{ \centerline{ \large Ilse Fischer} \vspace{0.3cm}
\centerline{Fakult\"at f\"ur Mathematik, Universit\"at Wien}
\centerline{Nordbergstrasse 15, A-1090 Wien, Austria}
\centerline{E-mail: {\tt Ilse.Fischer@univie.ac.at}} }

\begin{abstract} 
In the early 1980s, Mills, Robbins and Rumsey conjectured, and in 1996 Zeilberger proved a 
simple product formula for the number of $n \times n$ alternating sign matrices with 
a $1$ at the top of the $i$-th column. We give an alternative proof of this formula using 
our operator formula for the number of monotone triangles with prescribed bottom row. In addition, 
we provide the enumeration of certain $0$-$1$-($-1$) matrices generalizing alternating sign matrices.  
\end{abstract}

\maketitle

\section{Introduction}

An alternating sign matrix is a square matrix of $0$s, $1$s and $-1$s for which 
the sum of entries in each row and in each column is $1$ and the non-zero entries of each row and 
of each column alternate in sign. For instance, 
$$
\left(
\begin{array}{rrrrr}
0 & 0 & 1 & 0 & 0 \\
1 & 0 & -1 & 1& 0 \\
0 & 1 & 0 & -1 & 1 \\
0 & 0 & 1 & 0 & 0 \\
0 & 0 &  0 & 1 & 0
\end{array}
\right)
$$
is an alternating sign matrix. In \cite{mills2,mills} Mills, Robbins and 
Rumsey conjectured that there are 
\begin{equation}
\label{asmformula}
\prod_{j=1}^{n} \frac{(3j-2)!}{(n+j-1)!}
\end{equation}
$n \times n$ alternating sign matrices. This was first proved by Zeilberger~\cite{zeilberger}. 
Another, shorter, proof was given by Kuperberg~\cite{kuperberg} using the equivalence of alternating 
sign matrices with the six-vertex model for ``square ice'', which was earlier introduced in statistical 
mechanices. Zeilberger~\cite{zeilberger2} then used Kuperberg's observations to prove the following refinement generalizing 
\eqref{asmformula}.

\begin{theo}
\label{main}
The number of $n \times n$ 
alternating sign matrices where the unique $1$ in the first row is at 
the top of the $i$-th column is 
\begin{equation}
\label{refined}
\frac{(i)_{n-1}
(1+n-i)_{n-1}}{(n-1)!} \prod\limits_{j=1}^{n-1} \frac{(3j-2)!}{(n+j-1)!}.
\end{equation}
In this formula, $(a)_n = \prod\limits_{i=0}^n (a+i)$.
\end{theo}

The task of this paper is to give an alternative proof of Theorem~\ref{main}. As a byproduct, we 
obtain the  enumeration of certain objects generalizing alternating sign matrices.

\begin{theo}
\label{side}
The number of $n \times k$ matrices of $0$s, $1$s and $-1$s for which the non-zero entries 
of each row and each column alternate in sign and the sum in each row and each column is $1$, 
except for the columns $n, n+1,\ldots, k-1$, where we have sum $0$ and the first non-zero element is 
a $1$, is
$$
\prod_{j=1}^{n-1} \frac{(3j-2)!}{(n+j-1)!} \sum_{i=1}^n \frac{(i)_{n-1} (1+n-i)_{n-1}}{(n-1)!} 
\binom{i+k-n-1}{i-1}.
$$
\end{theo}

Our proofs are based on a 
formula \cite[Theorem~1]{fischer3} for the number of monotone triangles with 
given bottom row $(k_1,k_2,\ldots,k_n)$, which we have recently derived.
Monotone triangles with bottom row $(1,2,\ldots,n)$ are in bijection with $n \times n$ alternating 
sign matrices. The monotone triangle corresponding to a given alternating sign matrix
can be obtained as follows:
Replace every entry in the matrix by the sum of elements in the same column above, the 
entry itself included. In our running example we have
$$
\left(
\begin{array}{rrrrr}
0 & 0 & 1 & 0 & 0 \\
1 & 0 & 0 & 1 & 0 \\
1 & 1 & 0 & 0 & 1 \\
1 & 1 & 1 & 0 & 1 \\
1 & 1 &  1 & 1 & 1
\end{array}
\right).
$$
Row by row we record the columns that contain a $1$ and obtain the following triangular array.
$$
\begin{array}{ccccccccc}
  &  &  & & 3  &  &  &  &  \\
  &  &  &1 &   &  4&  &  &  \\
  &  &1  & &  2 &  & 5 &  &  \\
  &1  &  & 2&   & 3 &  &  5&  \\
1 & & 2 & & 3 & & 4 & & 5
\end{array}
$$
This is the monotone triangle corresponding to the alternating sign matrix above. Observe that 
it is weakly increasing in northeast direction and in southeast direction. Moreover, it is strictly 
increasing along rows. In general, a monotone triangle with $n$ rows is a triangular array 
$(a_{i,j})_{1 \le j \le i \le n}$ of integers such that $a_{i,j} \le a_{i-1,j} \le a_{i,j+1}$ 
and $a_{i,j} < a_{i,j+1}$ for all $i, j$. It is not too hard to see that 
monotone triangles with $n$ rows and $a_{n,j}=j$ 
are in bijection with $n \times n$ alternating sign matrices. Moreover, monotone triangles 
with $n-1$ rows and bottom row $(1,2,\ldots,i-1,i+1,\ldots,n)$ are in bijection with 
$n \times n$ alternating sign matrices with a $1$ at the bottom (or equivalently the top) of the $i$-th column. 
In \cite{fischer3} we gave the following operator formula for the number of monotone triangles with prescribed 
bottom row $(k_1,k_2,\ldots,k_n) \in \mathbb{Z}^n$.

\begin{theo}[\cite{fischer3}, Theorem~1]
\label{operatorformula}
The number of monotone triangles with $n$ rows and prescribed bottom row $(k_1, k_2, \ldots, k_n)$ 
is given by 
$$
\left( \prod_{1 \le p < q \le n} \left( \id + E_{k_p} \Delta_{k_q} \right) \right) 
\prod_{1 \le i < j \le n} \frac{k_j - k_i}{j-i}, 
$$
where $E_x$ denotes the shift operator, defined by 
$E_x \, p(x) = p(x+1)$, and $\Delta_x:= E_x - \id$ denotes the difference operator.
(In this formula, the product of operators is understood as the composition.)
\end{theo}

Thus, for instance, the number of monotone triangles with bottom row $(k_1,k_2,k_3)$
is 
\begin{multline*}
(\id + E_{k_1} \Delta_{k_2} ) (\id + E_{k_1} \Delta_{k_3}) (\id + E_{k_2} \Delta_{k_3})
\frac{1}{2} (k_2 - k_1) (k_3 - k_1) (k_3 - k_2) \\
= \frac{1}{2} ( -3 k_1 + k_1^2 + 2 k_1  k_2 - k_1^2 k_2 -
        2 k_2^2 + k_1  k_2^2 + 3 k_3 - 4 k_1  k_3  \\ +
        k_1^2 k_3 +   2 k_2 k_3  - k_2^2 k_3 + k_3^2 -
        k_1 k_3^2 + k_2 k_3^2).
\end{multline*}

We outline our proof of Theorem~\ref{main}. To prove the theorem using the formula in Theorem~\ref{operatorformula}
clearly means that we have to evaluate the formula at $(k_1,k_2,\ldots,k_n)=(1,2,\ldots,i-1,i+1,\ldots,n+1)$.
Let $A_{n,i}$ denote the number of $n \times n$ alternating 
sign matrices with a $1$ at the top of the $i$-th column and let $\alpha(n;k_1,\ldots,k_n)$ denote the 
number of monotone triangles with bottom row $(k_1,\ldots,k_n)$. Using the formula in Theorem~\ref{operatorformula},
we extend the interpretation of $\alpha(n;k_1,\ldots,k_n)$ to arbitrary $(k_1,\ldots,k_n) \in \mathbb{Z}^n$. 
In our proof, we first give a formula for $\alpha(n;1,2,\ldots,n-1,k)$ in terms of $A_{n,i}$. Next we show 
that $\alpha(n;1,2,\ldots,n-1,k)$ is an even polynomial in $k$ if $n$ is odd and an 
odd polynomial if $n$ is even. These two facts will then imply that $(A_{n,i})_{1 \le i \le n}$ is an 
eigenvector with respect to the eigenvalue $1$ of a certain matrix with binomial coefficients as entries. 
Finally we see that this determines $A_{n,i}$ up to a constant, which can easily be computed by induction with respect to $n$.

We believe that this second approach to prove the refined alternating sign matrix theorem not 
only provides us with a better understanding of known theorems, but will also enable us to obtain
new results in the field of plane partition and alternating sign matrix enumeration in the future.
A general strategy might be to derive analogous multivariate operators formulas for other
(triangular) arrays of integers, which correspond to certain classes of plane partitions 
and alternating sign matrices and which simplify to nice 
product formulas if we specialize the parameters in the right way. Hopefully these operator formulas can then 
also be used to derive these product formulas.

\section{A formula for $\alpha(n;1,2,\ldots,n-1,k)$}

We start by stating a fundamental recursion for $\alpha(n;k_1,\ldots,k_n)$. If we delete the last row of a monotone triangle 
with bottom row $(k_1,k_2,\ldots,k_n)$, we obtain 
a monotone triangle with $n-1$ rows and bottom row, say, $(l_1, l_2, \ldots, l_{n-1})$. By the definition of a monotone triangle, we have $k_1 \le l_1 \le k_2 \le l_2 \le \ldots \le k_{n-1} \le l_{n-1} \le k_n$ and $l_i \not= l_{i+1}$. Thus 
\begin{equation}
\label{rec}
\alpha(n;k_1,\ldots,k_n) = 
\sum_{(l_1,\ldots,l_{n-1}) \in \mathbb{Z}^{n-1}, \atop  k_1
\le l_1 \le k_2 \le \ldots \le k_{n-1} \le l_{n-1} \le
k_{n}, l_i \not= l_{i+1}} \alpha(n-1;l_1,\ldots,l_{n-1}).
\end{equation}

In the following lemma, we explain the action of operators, which
are symmetric polynomials in $E_{k_1}, E_{k_2}, \ldots, E_{k_n}$, 
on $\alpha(n;k_1,\ldots,k_n)$. It will be used twice in our proof of 
Theorem~\ref{main}.

\begin{lem} 
\label{sym-lem}
Let $P(X_1,\ldots,X_n)$ be a symmetric polynomial in $(X_1,\ldots,X_n)$ over $\mathbb{C}$. Then 
$$
P(E_{k_1},\ldots, E_{k_n}) \alpha(n;k_1,\ldots,k_n) = 
P(1,1,\ldots,1) \cdot \alpha(n;k_1,\ldots,k_n).
$$
\end{lem}

{\it Proof.} By Theorem~\ref{operatorformula} and the fact that shift operators with respect 
to different variables commute, it suffices to show that 
$$
P(E_{k_1},\ldots, E_{k_n}) 
\prod_{1 \le i < j \le n} \frac{k_j - k_i}{j-i} = 
P(1,1,\ldots,1) \cdot \prod_{1 \le i < j \le n} \frac{k_j - k_i}{j-i}.
$$
Let $(m_1,\ldots,m_n) \in \mathbb{Z}^n$ be with $m_i \ge 0$ for all $i$ and $m_i \not= 0$ for at 
least one $i$. It suffices to show that 
$$
\sum_{\pi \in {\mathcal S}_n} \Delta^{m_{\pi(1)}}_{k_1} \Delta^{m_{\pi(2)}}_{k_2} \dots 
\Delta^{m_{\pi(n)}}_{k_n}  \prod_{1 \le i < j \le n} \frac{k_j - k_i}{j-i} = 0. $$
By the Vandermonde determinant evaluation, we have 
$$
\prod_{1 \le i < j \le n} \frac{k_j - k_i}{j-i} = \det_{1 \le i, j \le n} \left( \binom{k_i}{j-1} \right).
$$
Therefore, it suffices to show that 
$$
\sum_{\pi, \sigma \in {\mathcal S}_n} \sgn \sigma \binom{k_1}{\sigma(1)-m_{\pi(1)}-1}
\binom{k_2}{\sigma(2)-m_{\pi(2)}-1} \dots \binom{k_n}{\sigma(n)-m_{\pi(n)}-1}=0.
$$
If, for fixed $\pi, \sigma \in {\mathcal S}_n$,  there exists an $i$ with $\sigma(i)-m_{\pi(i)}-1 < 0$ then the corresponding 
summand vanishes. We define a sign reversing involution on the set of non-zero summands. 
Fix $\pi, \sigma \in {\mathcal S}_n$ such that the summand corresponding to $\pi$ and $\sigma$ 
does not vanish. Consequently, $\{\sigma(1)-m_{\pi(1)}-1,\sigma(2)-m_{\pi(2)}-1,\ldots,\sigma(n)-m_{\pi(n)}-1\} \subseteq \{0,1,\ldots,n-1\}$ and since $(m_1,\ldots,m_n) \not= (0,\ldots,0)$, there are $i, j$, $1 \le i < j \le n$, 
with $\sigma(i)-m_{\pi(i)}-1=\sigma(j)-m_{\pi(j)}-1$. Among all pairs $(i,j)$ with this property, let $(i',j')$ be the pair, which is 
minimal with respect to the lexicographic order. Then the summand corresponding to 
$\pi \circ (i',j')$ and $\sigma \circ (i',j')$ is the negativ of the summand corresponding to $\pi$ and $\sigma$. \qed

\medskip

Let 
$$
e_p(X_1,\ldots, X_n) = \sum_{1 \le i_1 < i_2 < \ldots i_p \le n} X_{i_1} X_{i_2} \dots X_{i_p}
$$
denote the $p$-th elementary symmetric function. Lemma~\ref{sym-lem}  will be used to 
deduce a formula for 
\begin{equation}
\label{A}
\left. e_{p-j}(E_{k_1},E_{k_2},\dots,E_{k_{n-1}}) \alpha(n;k_1,\ldots,k_n) 
\right|_{(k_1,\ldots,k_n)=(1,2,\ldots,n-1,n+j)}
\end{equation}
in terms of $A_{n,i}$  if $0 \le j \le p \le n-1$ (Lemma~\ref{last}). If we specialize $p=j$ in this 
identity we obtain the 
desired formula for $\alpha(n;1,2,\ldots,n-1,k)$, which we have mentioned in the outline of our proof. 
The formula for \eqref{A} will be shown by 
induction with respect to $j$. In the following lemma, we deal with the initial case of the induction.

\begin{lem}
\label{jzero}
Let $0 \le p \le n-1$. Then we have
$$
\left. e_p(E_{k_1},E_{k_2},\dots,E_{k_{n-1}}) \, \alpha(n;k_1,\ldots,k_n) \right|_{(k_1,\ldots,k_n)=(1,2,\ldots,n)} 
= \sum_{i=1}^n \binom{n-i}{p} A_{n,i}.
$$
\end{lem}

{\it Proof.} First observe that for $1 \le i_1 < i_2 < \ldots i_p \le n-1$ we have 
\begin{multline*}
\left. E_{k_{i_1}} E_{k_{i_2}} \ldots E_{k_{i_p}} \alpha(n;k_1,\ldots,k_n) \right|_{(k_1,\ldots,k_n)=(1,2,\ldots,n)} \\
= \alpha(n;1,2,\ldots, i_{1}-1, i_1+1, i_1+1, \ldots, i_2-1, i_2+1, i_2+1, \ldots, i_p-1, i_p+1, i_p+1, \ldots, n) \\
= \sum_{1 \le j_1 \le 2 \le j_2 \le \ldots \le i_{1}-1 \le j_{i_{1}-1} \le i_{1}+1 \atop j_1 < j_2 < \ldots < j_{i_{1}-1}} 
\alpha(n-1;j_1,j_2,\ldots, j_{i_{1}-1},i_{1}+1,i_{2}+1,\ldots,n) \\
= \alpha(n;1,2,\ldots, i_{1}-1, i_1+1, i_1+1, i_1+2, i_1+3, \ldots, n),
\end{multline*}
where the second equality follow from \eqref{rec}.
In particular, we see that 
$$\left. E_{k_{i_1}} E_{k_{i_2}} \ldots E_{k_{i_p}} \alpha(n;k_1,\ldots,k_n) \right|_{(k_1,\ldots,k_n)=(1,2,\ldots,n)}$$ does not depend on $i_2,\ldots, i_p$. 
Consequently, the left-hand-side in the statement of the lemma is equal to 
\begin{equation}
\label{1}
\sum_{j=1}^{n-1} \binom{n-1-j}{p-1} \alpha(n;1,2,\ldots,j-1,j+1,j+1,j+2,\ldots,n).
\end{equation}
Next observe that by \eqref{rec}
\begin{multline*}
\alpha(n;1,2,\ldots,j-1,j+1,j+1,j+2,\ldots,n) \\
= \sum_{i=1}^{j} \alpha(n-1;1,2,\ldots,i-1,i+1,\ldots,n) = \sum_{i=1}^{j} A_{n,i}.
\end{multline*}
Thus, \eqref{1} is equal to 
$$
\sum_{j=1}^{n-1} \binom{n-1-j}{p-1} \sum_{i=1}^{j} A_{n,i} 
= \sum_{i=1}^n \sum_{j=i}^n \binom{n-1-j}{p-1} A_{n,i} - \binom{-1}{p-1} \sum_{i=1}^{n} A_{n,i}.
$$
We complete the proof by using the following summation formula
$$
\sum_{j=a}^b \binom{x+j}{n} = \sum_{j=a}^b \left( \binom{x+j+1}{n+1} - \binom{x+j}{n+1} \right) = \binom{x+b+1}{n+1} - \binom{x+a}{n+1}.
$$
\qed

\begin{lem}
\label{last}
Let $0 \le j \le p \le n-1$. Then
\begin{multline*}
\left. e_{p-j}(E_{k_1},E_{k_2},\dots,E_{k_{n-1}}) \alpha(n;k_1,\ldots,k_n) \right|_{(k_1,\ldots,k_n)=(1,2,\ldots,n-1,n+j)} \\
= (-1)^j \sum_{i=1}^n A_{n,i} \left( \binom{n-i}{p} + \sum_{l=0}^{j-1} \binom{n}{p-l} \binom{i+l-1}{l} (-1)^{l-1} \right)
\end{multline*}
For $p=j$ this simplifies to 
$$
\alpha(n;1,\ldots,n-1,n+j) = \sum_{i=1}^n A_{n,i} \binom{i+j-1}{i-1}.
$$ 
\end{lem}

{\it Proof.} First we show how the first formula implies the second. For this 
purpose, we have to consider 
$$
(-1)^j \left( \binom{n-i}{j} + \sum_{l=0}^{j-1} \binom{n}{j-l} \binom{i+l-1}{l} (-1)^{l-1} \right).
$$
Using $\binom{i+l-1}{l} = \binom{-i}{l} (-1)^l$, we see that this is equal to 
$$
(-1)^j \left( \binom{n-i}{j} + (-1)^j \binom{i+j-1}{j} - \sum_{l=0}^j \binom{n}{j-l} \binom{-i}{l} \right).
$$
We apply the Chu-Vandermonde summation \cite[p. 169, (5.22)]{knuth}, in order to see that this is simplifies to $\binom{i+j-1}{j}=\binom{i+j-1}{i-1}$.

Now we consider the first formula. By Lemma~\ref{sym-lem}, we have 
\begin{multline*}
\left. e_{p-j} (E_{k_1},\ldots,E_{k_n}) \alpha(n;k_1,\ldots,k_n) \right|_{(k_1,\ldots,k_n)=(1,2,\ldots,n-1,n+j)} \\ = 
\binom{n}{p-j} \alpha(n;1,2,\ldots,n-1,n+j).
\end{multline*}
On the other hand, we have 
\begin{multline*}
\left. e_{p-j} (E_{k_1},\ldots,E_{k_n}) \alpha(n;k_1,\ldots,k_n) \right|_{(k_1,\ldots,k_n)=(1,2,\ldots,n-1,n+j)} \\
= \left. e_{p-j-1} (E_{k_1},\ldots,E_{k_{n-1}}) \alpha(n;k_1,\ldots,k_n)\right|_{(k_1,\ldots,k_n)=(1,2,\ldots,n-1,n+j+1)}\\
+ \left.  e_{p-j} (E_{k_1},\ldots,E_{k_{n-1}}) \alpha(n;k_1,\ldots,k_n)\right|_{(k_1,\ldots,k_n)=(1,2,\ldots,n-1,n+j)}.
\end{multline*}
This implies the recursion 
\begin{multline*}
\left. e_{p-j-1} (E_{k_1},\ldots,E_{k_{n-1}}) \alpha(n;k_1,\ldots,k_n)\right|_{(k_1,\ldots,k_n)=(1,2,\ldots,n-1,n+j+1)} \\
= \binom{n}{p-j} \alpha(n;1,\ldots,n-1,n+j)  \\ - 
\left.  e_{p-j} (E_{k_1},\ldots,E_{k_{n-1}}) \alpha(n;k_1,\ldots,k_n)\right|_{(k_1,\ldots,k_n)=(1,2,\ldots,n-1,n+j)}.
\end{multline*}
Now we can prove the first formula in the lemma by induction with respect to $j$. The case $j=0$ was dealt with 
in Lemma~\ref{jzero}. \qed

\medskip

By Theorem~\ref{operatorformula}, $\alpha(n;1,2,\ldots,n-1,n+j)$ is a polynomial in $j$ of degree no greater than $n-1$. 
By Lemma~\ref{last}, it coincides with $\sum\limits_{i=1}^n A_{n,i} \binom{i+j-1}{i-1}$ for $j=0,1,\ldots, n-1$
and, since $\sum\limits_{i=1}^n A_{n,i} \binom{i+j-1}{i-1}$ is a polynomial in $j$ of degree no greater than $n-1$ as well, the two polynomials
must be equal. This constitutes the following.

\begin{lem} 
\label{k}
Let $n \ge 1$ and $k \in \mathbb{Z}$. Then we have 
$$
\alpha(n;1,2,\ldots,n-1,k) = \sum_{i=1}^n A_{n,i} \binom{i+k-n-1}{i-1}.
$$
\end{lem}

\section{The symmetry of $k \to \alpha(n;1,2,\ldots,n-1,k)$}

In the following lemma we prove a transformation formula for $\alpha(n;k_1,\ldots,k_n)$, which  
implies as a special case that $k \to \alpha(n;1,2,\ldots,n-1,k)$ is even if $n$ is odd and 
odd otherwise.

\begin{lem}
\label{shift}
Let $n \ge 1$. Then we have 
$$
\alpha(n;k_1,\ldots,k_n) = (-1)^{n-1} \alpha(n;k_2,\ldots,k_n, k_1 - n).
$$
\end{lem}

{\it Proof.} By Theorem~\ref{operatorformula}
\begin{multline*}
(-1)^{n-1} \alpha(n;k_2,\ldots,k_n,k_1-n) \\
= (-1)^{n-1} \prod_{2 \le p < q \le n} (\id + E_{k_p} \Delta_{k_q}) 
\prod_{p=2}^n (\id + E_{k_p} \Delta_{k_1}) \prod_{2 \le i < j \le n} \frac{k_j - k_i}{j-i} 
\prod_{i=2}^n \frac{k_1 - k_i - n}{i-1} \\
= \prod_{2 \le p < q \le n} (\id + E_{k_p} \Delta_{k_q}) 
\prod_{p=2}^n (\id + E_{k_p} \Delta_{k_1}) \prod_{2 \le i < j \le n} \frac{k_j - k_i}{j-i} 
\prod_{i=2}^n \frac{k_i + n - k_1}{i-1} \\
= \prod_{2 \le p < q \le n} (\id + E_{k_p} \Delta_{k_q}) 
\prod_{p=2}^n (\id + E_{k_p} \Delta_{k_1}) E^{-n}_{k_1} \prod_{1 \le i < j \le n} \frac{k_j - k_i}{j-i}.
\end{multline*}
Thus, we have to show that 
\begin{multline*}
\prod_{2 \le p < q \le n} (\id + E_{k_p} \Delta_{k_q}) 
\left( \prod_{q=2}^{n} \left( \id + E_{k_1} \Delta_{k_q} \right) -
E^{-n}_{k_1} \prod_{q=2}^n (\id + E_{k_q} \Delta_{k_1}) \right) 
\prod_{1 \le i < j \le n} \frac{k_j - k_i}{j-i} = 0.
\end{multline*}
Clearly, it suffices to prove that
$$
\left( E_{k_1}^n \prod_{q=2}^{n} \left( \id + E_{k_1} \Delta_{k_q} \right) -
 \prod_{q=2}^n (\id + E_{k_q} \Delta_{k_1}) \right) 
\prod_{1 \le i < j \le n} \frac{k_j - k_i}{j-i} = 0.
$$
We replace $\Delta_{k_i}$ by $X_i$ and, accordingly, 
$E_{k_i}$ by $(X_i + 1)$ in the operator in this expression and obtain 
\begin{equation}
\label{poly}
(X_1+1)^n \prod_{q=2}^n (1+ (X_1 +1) X_q) - \prod_{q=2}^{n} (1 + (X_q + 1) X_1).
\end{equation}
By the proof of Lemma~\ref{sym-lem}, the assertion follows if we show that this 
polynomial 
is in the ideal, which is generated by the symmetric polynomials in $X_1, X_2, \ldots, X_n$ 
without constant term. Observe that \eqref{poly} is equal to 
\begin{multline}
\label{calc1}
\sum_{j=0}^{n-1} (X_1+1)^{n+j} e_j(X_2,\ldots,X_n) - 
\sum_{j=0}^{n-1} X_1^j \, e_j(X_2+1,X_3+1,\ldots,X_n+1) \\
= 
\sum_{j=0}^{n-1} (X_1+1)^{n+j} e_j(X_2,\ldots,X_n) - 
\sum_{j=0}^{n-1} X_1^j \sum_{i=0}^j \binom{n-i-1}{j-i} 
e_i(X_2,\ldots,X_n) \\
=
\sum_{j=0}^{n-1} \left( (X_1+1)^{n+j} - X_1^j (X_1+1)^{n-j-1} \right) e_j(X_2,\ldots,X_n).
\end{multline}
We recursively define a sequence $(q_j(X))_{j \ge 0}$ of Laurent polynomials. Let
$q_0(X)=0$ and 
\begin{equation}
\label{laurent-rec}
q_{j+1}(X) = (X+1)^{2j+1} - X^j - q_j(X) - q_j(X) \, X^{-1}.
\end{equation}
We want to show that this is in fact a sequence of polynomials having a 
zero at $X=0$. For this purpose, we consider 
$$
Q(X,Y) := \sum_{j \ge 0} q_j(X) Y^j.
$$
Using \eqref{laurent-rec} and the initial condition, we obtain the following
$$
Q(X,Y) = \frac{X \, Y}{(1-X \, Y) ( 1 - (X+1)^2 Y)},
$$
which immediately implies the assertion. We set 
$$
p_{j}(X) = q_{j}(X) (X+1)^{n-j} X^{-1}
$$
and observe that, for all $j$ with $j \le n$, $p_{j}(X)$ is a polynomial in $X$.
The recursion \eqref{laurent-rec} clearly implies 
$$
p_{j+1}(X) X =  (X+1)^{n+j} - X^j (X+1)^{n-j-1} - p_j(X).
$$
Thus, \eqref{calc1} is equal to 
\begin{equation*}
\sum_{j=0}^{n-1} \left( p_j(X_1) + p_{j+1}(X_1) X_1 \right) e_j(X_2,\ldots,X_n) \\
= \sum_{j=0}^n p_j(X_1) e_j(X_1,\ldots,X_n).
\end{equation*}
Since $p_0(X)=0$, this expression is in the ideal generated by 
the symmetric polynomials in $(X_1,\ldots, X_n)$ without constant term and the assertion 
of the lemma is proved. \qed

\medskip

If $(a_{i,j})_{1 \le j \le i \le n}$ is a monotone triangle 
with bottom row $(k_1,\ldots,k_n)$  then $(-a_{i,n+1-j})_{1 \le j \le i \le n}$
is a monotone triangle with bottom row $(-k_n,\ldots,-k_1)$.
This implies the following identity.
\begin{equation}
\label{negativ}
\alpha(n;k_1,k_2,\ldots,k_n) = \alpha(n;-k_n,-k_{n-1},\ldots,-k_1)
\end{equation}
Similarly, it is easy to see that 
\begin{equation}
\label{constant}
\alpha(n;k_1,k_2,\ldots,k_n) = \alpha(n;k_1+c,k_2+c,\ldots,k_n+c)
\end{equation}
for every integer constant $c$. Therefore
\begin{multline*}
\alpha(n;1,2,\ldots,n-1,k) 
= \alpha(n;-k,-n+1,-n+2,\ldots,-1) \\
= (-1)^{n-1} \alpha(n;-n+1,-n+2,\ldots,-1,-k-n) = (-1)^{n-1} \alpha(n;1,2,\ldots,n-1,-k),
\end{multline*}
where the first equality follows from \eqref{negativ}, the second from Lemma~\ref{shift} and 
the third from \eqref{constant} with $c=n$. This, together with Lemma~\ref{k}, implies the following 
identity
$$
\sum_{i=1}^n A_{n,i} \binom{i+k-n-1}{i-1} = (-1)^{n-1} \sum_{i=1}^n A_{n,i} \binom{i-k-n-1}{i-1}
$$
for all integers $k$.
In this identity, we replace $\binom{i-k-n-1}{i-1}$ by $(-1)^{i-1} \binom{k+n-1}{i-1}$ and 
$\binom{i+k-n-1}{i-1}$ by 
$$
\sum_{j=1}^{n} \binom{i-2n}{i-j} \binom{k+n-1}{j-1},
$$
which is possible by the Chu-Vandermonde summation \cite[p. 169, (5.22)]{knuth} if $i \le n$. We interchange the role of $i$ and $j$ on the 
left-hand-side and obtain
$$
\sum_{i=1}^n \sum_{j=1}^n A_{n,j} \binom{j-2n}{j-i} \binom{k+n-1}{i-1} = 
\sum_{i=1}^n A_{n,i} (-1)^{n+i} \binom{k+n-1}{i-1}
$$
for all integers $k$.
Since $(\binom{k+n-1}{i-1})_{i \ge 1}$ is a basis of $\mathbb{C}[k]$ as a vectorspace over 
$\mathbb{C}$, this implies 
$$
\sum_{j=1}^n A_{n,j} \binom{j-2n}{j-i}  = A_{n,i} (-1)^{n+i}
$$
for $i=1,2,\ldots,n$.
We replace $\binom{j-2n}{j-i}$ by $(-1)^{j-i} \binom{2n-i-1}{j-i}$. Moreover 
we replace $j$ by $n+1-j$ and use the fact that $A_{n,j}=A_{n,n+1-j}$ in order 
to obtain 
\begin{equation}
\label{eigenvector}
\sum_{j=1}^n A_{n,j} (-1)^{j+1} \binom{2n-i-1}{n-i-j+1} = A_{n,i}.
\end{equation}
Phrased differently, $(A_{n,1},A_{n,2},\ldots,A_{n,n})$ is an eigenvector of 
$((-1)^{j+1} \binom{2n-i-1}{n-i-j+1})_{1 \le i,j \le n}$ with respect to the 
eigenvalue $1$. In the following section we see that this determines 
$(A_{n,1},A_{n,2},\ldots,A_{n,n})$ up to a multiplicative constant, which 
we are able to compute.

\section{The eigenspace of $ \left( (-1)^{j+1} \binom{2 n - i -1}{n-i-j+1} \right)_{1 \le i , j \le n}$ with respect to the eigenwert $1$ is one dimensional.}

Since $(A_{n,1},\ldots,A_{n,n})$ is an eigenvector of $ \left( (-1)^{j+1} \binom{2 n - i -1}{n-i-j+1} \right)_{1 \le i , j \le n}$, it suffices to show that the dimension of the eigenspace with respect to $1$ is no greater than $1$.
Thus we have to show that the rank of
$$
\left( (-1)^{j} \binom{2 n - i -1}{n-i-j+1} + \delta_{i,j} \right)_{1 \le i , j \le n}
$$
is at least $n-1$.  It suffices to show that 
$$
\det_{2 \le i,j \le n} \left( (-1)^{j} \binom{2 n - i -1}{n-i-j+1} + \delta_{i,j} \right) \not = 0.
$$
We shift $i$ and $j$ by one in this determinant and obtain
\begin{equation}
\label{det}
\det_{1 \le i,j \le n-1} \left( (-1)^{j+1} \binom{2 n - i -2}{n-i-j-1} + \delta_{i,j} \right).
\end{equation}
Let $B_n$ denote the matrix underlying the determinant.
We define $R_n = \left( \binom{n+j-i-1}{j-i} \right)_{1 \le i , j \le n-1}$.  Observe that 
$R_n^{-1} = \left( (-1)^{i+j} \binom{n}{j-i} \right)_{1 \le i , j \le n-1}$. 
Moreover, we have 
$
R_n^{-1} B_n R_n = B^{*}_n + I_{n-1},
$
where $I_{n-1}$ denotes the $(n-1) \times (n-1)$ identity matrix and $B^{*}_n$
is the $(n-1) \times (n-1)$ matrix with $\binom{i+j}{j-1}$ as entry in the $i$-th row 
and $j$-th column except for the last row, where we have all zeros. (This transformation 
is due to Mills, Robbins and Rumsey \cite{mills2}.) Thus the 
determinant in \eqref{det} is equal to 
$$
\det (B^{*}_n + I_{n-1}) = \det_{1 \le i , j \le n-2} \left( \binom{i+j}{j-1} + I_{n-2} \right),
$$  
where we have expanded $B^{*}_n + I_{n-1}$ with respect to the last row. Andrews \cite{andrews}
has shown that this determinant gives the number of desecending plane partitions with no part greater than 
$n-1$ and, therefore, the determinant does not vanish. Recently, Krattenthaler \cite{kratt} showed
that descending plane partitions can be geometrically realized as cyclically symmetric rhombus tilings of a certain hexagon of which a centrally located equilateral triangle of side length 2 has been removed.

In order to complete our proof, we have to show that $(B_{n,i})_{1 \le i \le n}$ with 
$$
B_{n,i}= \frac{ (i)_{n-1} (1+n-i)_{n-1}}{(n-1)!} \prod_{k=1}^{n-1} \frac{(3k-2)!}{(n+k-1)!} 
$$
is an eigenvector of $ \left( (-1)^{j+1} \binom{2 n - i -1}{n-i-j+1} \right)_{1 \le i , j \le n}$ 
with respect to the eigenwert $1$, i.e. we have to show that
\begin{multline*}
\sum_{j=1}^{n} (-1)^{j+1} \binom{2n-i-1}{n-i-j+1} \frac{(j)_{n-1} (1+n-j)_{n-1} }{(n-1)!} 
\prod_{k=1}^{n-1} \frac{(3k-2)!}{(n+k-1)!} = \\
\frac{(i)_{n-1} (1+n-i)_{n-1}}{(n-1)!} \prod_{k=1}^{n-1} \frac{(3k-2)!}{(n+k-1)!}.
\end{multline*}
This is equivalent to showing that 
\begin{equation}
\label{identity}
\sum_{j=1}^n (-1)^{j+1} \binom{2n-i-1}{n-j-i+1} \binom{n+j-2}{n-1} \binom{2n-j-1}{n-1} 
= 
\binom{n+i-2}{n-1} \binom{2n-i-1}{n-1}.
\end{equation}
Observe that the left hand side of this identity is equal to 
\begin{multline*}
\binom{2n-i-1}{n-1} \binom{n+i-2}{n-1} \binom{n-1}{i-1}^{-1} 
\sum_{j=1}^{n} (-1)^{j+1} \binom{2n-j-1}{n-j-i+1} \binom{n-1}{j-1} = \\ 
\binom{2n-i-1}{n-1} \binom{n+i-2}{n-1} \binom{n-1}{i-1}^{-1}
\sum_{j=1}^{n} (-1)^{n+i} \binom{-n-i+1}{n-j-i+1} \binom{n-1}{j-1} = \\
\binom{2n-i-1}{n-1} \binom{n+i-2}{n-1} \binom{n-1}{i-1}^{-1} 
(-1)^{n+i} \binom{-i}{n-i} = 
\binom{2n-i-1}{n-1} \binom{n+i-2}{n-1},
\end{multline*} 
where the first and third equality follows from 
$\binom{a}{b}= (-1)^b \binom{b-a-1}{b}$ and the 
second equality follows from the Chu-Vandermonde
identity; see \cite[p. 169, (5.26)]{knuth}. Consequently, $A_{n,i} = C_n \cdot B_{n,i}$ for 
$C_n \in \mathbb{Q}$.
It is easy to check that $C_1=1$. Observe that $\sum\limits_{i=1}^{n-1} A_{n-1,i} = A_{n,1}$, 
since $(n-1) \times (n-1)$ alternating sign matrices are bijectively related to 
$n \times n$ alternating sign matrices with a $1$ at the top of the first column.
Moreover, observe that we also have $\sum\limits_{i=1}^{n-1} B_{n-1,i} = B_{n,1}$, since 
\begin{multline*}
\sum_{i=1}^{n-1} B_{n-1,i} = 
\prod_{k=1}^{n-2} \frac{(3k-2)!}{(n+k-2)!} \sum_{i=1}^{n-1} \frac{(i)_{n-2} (n-i)_{n-2}}{(n-2)!}  = \\
\prod_{k=1}^{n-2} \frac{(3k-2)!}{(n+k-2)!} (n-2)! \sum_{i=1}^{n-1} \binom{i+n-3}{n-2} \binom{2n-i-3}{n-2} = \\
\prod_{k=1}^{n-2} \frac{(3k-2)!}{(n+k-2)!} (n-2)! \binom{3n-5}{2n-3} =
B_{n,1},
\end{multline*}  
where the second equality follows from an identity, which is equivalent to the  Chu-Vandermonde
identity; see \cite[p. 169, (5.26)]{knuth}. Therefore, by induction with respect to $n$, we have 
$C_n=1$ for all $n$. This completes our proof of the refined alternating sign matrix theorem. Theorem~\ref{side}
follows if we combine Theorem~\ref{main} and Lemma~\ref{k}, since a careful analysis of the bijection 
between alternating sign matrices and monotone triangles shows that $\alpha(n;1,2,\ldots,n-1,k)$
is the number of objects described in the statement of the theorem.

\end{document}